\documentclass[12pt,leqno]{article}
\title{On the class number of cyclic extensions $K/\Q$. Draft V5.0}
\author{Roland Qu\^eme}

\usepackage{amsmath,amsthm,amstext,amsbsy,amsfonts}
\newtheorem{thm}{Theorem}[section]
\newtheorem{cor}[thm]{Corollary}

\newtheorem{lem}[thm]{Lemma}
\newcommand{\N}{\mathbb{N}}

\newcommand{\C}{\mathbb{C}}
\newcommand{\Q}{\mathbb{Q}}
\newcommand{\Z}{\mathbb{Z}}
\newcommand{\s}{\mathbf{s}}
\newcommand{\modu}{\ \mbox{mod}\ }
\newcommand{\be}{\begin{equation}}
\newcommand{\ee}{\end{equation}}
\date{May 16 1999}
\begin{document}
\abstract

This draft is a contribution to the description of the prime factors of the class number of a cyclic extension $K/\Q$ :
\begin{itemize}
\item
An example  of the results obtained in the first approach is  : 
Let $K/\Q$ be a cyclic extension with $[K:\Q]=g$.
Let $g=2^\gamma\prod_{j=1}^m g_j^{\gamma_j}$ be the factorization of $g$ in primes. Suppose that $0\leq \gamma\leq 1$.
Let $G$ be the class group of $K/\Q$.
Let $h$ be the class number of $K/\Q$.
Let $h=\prod_{i=1}^r h_i^{\nu_i}$ be the factorization of $h$ in primes.
Let $\sigma$ be any  $\Q-$isomorphism generating $K/\Q$.
If for one prime factor $h_i$ of $h$, the $h_i$-component $G(h_i)$ of the class group $G$ of $K$  is cyclic, then  
\begin{enumerate}
\item
else $g\equiv 0 \mod h_i$,
\item
else $g\not\equiv 0 \modu h_i$ and 
\begin{enumerate}
\item
else $h_i\equiv 1 \modu g_j$ for a prime odd divisor $g_j$ of $g$, $1\leq j \leq m$,
\item
else  $g\equiv 0 \modu 2$ and $h_i$ divides the class number $h(L_2)$ of the quadratic extension  $L_2,\quad \Q\subset L_2\subset K$. In that case for all ideal $\mathbf s $ generating the cyclic subgroup $\Gamma(h_i)$ of order $h_i$ of the class group $G$, the ideal $\mathbf s\times\sigma(\mathbf s)$ is principal. 
\end{enumerate}
\end{enumerate}
\item An example of result obtained with a second approach is :
Let $p\in\N,\quad p>3$ be a prime.
Let $a,b\in\Z-\{0\},\quad (a+b)\not=0,\quad \frac{(a^p+b^p)}{(a+b)}>1$.
Suppose that $a,b$ are coprime with $ab(a+b)\not\equiv 0\modu p$. 
Let $n=p\frac{(a^p+b^p)}{(a+b)}$. Let $\Q(\xi_n)$ be the $n-$cyclotomic field.
Let $h_n$ be the class number of $\Q(\xi_n)$.
Then $h_n\equiv 0 \modu p$.
\end{itemize}
The proofs are {\bf stricly elementary}.
We give  several verifications of results obtained from the tables in  Washington\cite{was}, Schoof \cite{sch}, Masley \cite{mas}, Girstmair \cite{gir}, Jeannin \cite{jea} and the tables of cubic totally real number fields of the $ftp$ server
 {\it megrez.math.u-bordeaux.fr}.
\endabstract
\maketitle
%
\section{On class number of cyclic extensions  $K/\Q$.}\label{s22041}
\subsection{Some definitions and notations}\label{ss22041}
\begin{itemize}
\item
Let $K/\Q$ be a cyclic extension with $[K:\Q]=g$.
\item
Let $B$ be the ring of integers of $K$.
\item 
Let $g=2^{\gamma}\prod_{j=1}^s g_i^{\gamma_j}$ be the factorization of $g$ in primes. Suppose, moreover that $0\leq \gamma\leq 1$
\item 
Let $\sigma$ be a $\Q-$isomorphism which generates $Gal(K/\Q)$.
\item
Let $G= Cl(K)$ be the class group of $K/\Q$.
\item 
For an integral ideal $\mathbf s$, let us denote $Cl(\mathbf s)$, the class of $\mathbf s$ in the class group $G$.
\item
Let $h$ be the class number of $K/\Q$. 
\item
Let $h=\prod_{i=1}^r h_i^{\nu_i}$ be the factorization of $h$ in primes in $\N$.
\item 
Let us suppose that, for the prime $h_i$ studied,  the $h_i-$component $G(h_i)$ of the class-group $G$ of $K/\Q$  is cyclic. In the sequel we study the properties of this prime factor $h_i$ of $h$. In that case, there is one and only one subgroup $\Gamma(h_i)$ of $G(h_i)$ with order $h_i$. Note that if $h\not \equiv 0 \modu h_i^2$, then $G(h_i)$ is of order $h_i$, hence cyclic.
\item
When $g$ has at least two different odd prime divisors, let  $d,\delta\in \N, \quad 2 < d \leq g, 
\quad d\times \delta=g,\quad (d,\delta)=1$. In that case, let $L_d$ be the intermediate field, $\Q\subset L_d\subset K$ with $[L_d:\Q]=d$. Let $B_d$ be the ring of integers of $L_d$. The extension $L_d/\Q$ is cyclic of order $d$. The Galois group $Gal(K/\Q)=C_1\oplus C_2$ with $Card(C_1)=d$ and $Card(C_2)=\delta$ and $Gal(K/L_d)$ is also cyclic.
\end{itemize}
%
\subsection{Factorization of the class number $h$ of a cyclic extension $K/\Q$.}\label{ss22042}
\begin{lem}\label{l14051} Let $K$ be a number field of degree $n$, signature $(r_1,r_2)$ and absolute discriminant $D$. Let $B$ be the Minkoswki bound of $K$, given by $B=(\frac{4}{\pi})^{r_2}\times \frac{n!}{n^n}\times\sqrt{D}$. Let $h$ be the class number of $K$.  If $B>17$, then  $h<(B\times(2Log(B))^n$.
\begin{proof}
See Qu\^eme, \cite{que} Theorem 3.3 p 19-04.
\end{proof}
\end{lem}
\begin{thm}\label{t22041}
Let $K/\Q$ be a cyclic extension with $[K:\Q]=g$.
Let $g=2^\gamma\prod_{j=1}^m g_j^{\gamma_j}$ be the factorization of $g$ in primes. Suppose that $0\leq \gamma\leq 1$.
Let $G$ be the class group of $K/\Q$.
Let $h$ be the class number of $K/\Q$.
Let $h=\prod_{i=1}^r h_i^{\nu_i}$ be the factorization of $h$ in primes.
Let $\sigma$ be any  $\Q-$isomorphism generating $K/\Q$.
If for one prime factor $h_i$ of $h$, the $h_i$-component $G(h_i)$ of the class group $G$ of $K$  is cyclic, then  
\begin{enumerate}
\item
else $g\equiv 0 \mod h_i$,
\item
else $g\not\equiv 0 \modu h_i$ and 
\begin{enumerate}
\item
else $h_i\equiv 1 \modu g_j$ for a prime odd divisor $g_j$ of $g$, $1\leq j \leq m$,
\item
else  $g\equiv 0 \modu 2$ and $h_i$ divides the class number $h(L_2)$ of the quadratic field $L_2, \quad \Q\subset L_2\subset K$. In that case for all ideal $\mathbf s $ generating the cyclic subgroup $\Gamma(h_i)$ of order $h_i$ of the class group $G$, the ideal $\mathbf s\times\sigma(\mathbf s)$ is principal. 
\end{enumerate}
\end{enumerate}
\begin{proof}$ $
\begin{itemize}
\item
The subgroup $G(h_i)$ of $G$ is from hypothesis a cyclic group. 
It has one and only one subgroup $\Gamma(h_i)$ of order $h_i$. 
Let $\s$ be an integral ideal of $K$ not principal and  such that 
$\s^{h_i}$ is principal: note that such an ideal exists. 
We have $<Cl(\s)>=\Gamma(h_i)$. 
From $\s^{h_i}=B$, we get $\sigma(\s^{h_i})=B$ and hence $<Cl(\sigma(\s))>=\Gamma(h_i)$. 
Therefore, there exists $\alpha\in\N, \quad 0\leq \alpha\leq h_i-1$ such that $Cl(\sigma(\s))=Cl(\s^\alpha)$.
In the same way, we have $Cl((\sigma^2(\s))=Cl(\sigma\circ\sigma(\s))
=Cl(\sigma(\sigma(\s))
=Cl(\sigma(\s^\alpha))$ 
because 
$\frac{\sigma(\s)}{\s^\alpha}$ is a principal ideal, so $\sigma(\frac{\sigma(\s)}{\s^\alpha})$ is principal and so $Cl(\sigma^2(\s))=Cl(\sigma(\s^\alpha))$.
But $Cl(\sigma(\s^\alpha))
=\alpha\times Cl(\sigma(\s))
=\alpha\times Cl(\s^\alpha)
=Cl(\s^{\alpha^2})$ 
and finally $Cl(\sigma^2(\s))=Cl(\s^{\alpha^2})$.
Then, we have for $j=1,\dots,g-1$ the relation $Cl(\sigma^i(\s))=Cl(\s^{\alpha^i})$.
\item 
In an other part $N_{K/\Q}(\s)B=\prod_{i=0}^{g-1} \sigma^i(\s)$ is a principal ideal of $B$. Therefore 
$\s^{1+\alpha+\alpha^2+\dots+\alpha^{g-1}}$ is a principal ideal for some 
$\alpha\in\N,\quad 1\leq \alpha \leq h_i-1$.

\begin{itemize}
\item First case: $\alpha=1$

Then $\s^{g}$ is  a principal ideal  and $g\equiv 0 \modu h_i$.
\item Second case : $g\not\equiv 0 \modu 2$ and $\alpha\not=1$.
 
We have $\frac{\alpha^g-1}{\alpha-1}\equiv 0 \modu h_i$. From classical properties of Vandiver form $\frac{\alpha^g-1}{\alpha-1}$, we deduce immediatly that there exists at least one divisor $g_j$ of $g$ such that $h_i\equiv 1 \modu g_j$.
\item Third case : $2\|g$ and $\alpha\not\in \{1, h_i-1\}$. 

Then $\alpha^2-1\not\equiv 0 \modu h_i$, then
same result as second case with the form $\frac{\alpha^g-1}{\alpha^2-1}$.
\item Fourth case : $2\|g$ and $\alpha=h_i-1$. 

\begin{itemize}
\item 
Show at first that the ideal $\s\times \sigma(\s)$ is principal:
Then $Cl(\sigma(\mathbf s))=Cl(\mathbf s^{h_i-1})=Cl(\mathbf s^{-1})$ and thus $\mathbf s\times\sigma(\mathbf s)$ is a principal ideal for all $\Q-$isomorphisms which generate the class group $G$ of $K/\Q$.
\item  Let us consider the ideal 
$\mathbf b=N_{K/L_2}(\s)B=\prod_{i=0}^{g/2-1} \sigma^{2i}(\s)$. 
In the field $K$, we have from above $Cl(\s)=Cl(\sigma^2(\s)=\dots=Cl(\sigma^{g/2-1}(s))$. Therefore,  we get $Cl(\mathbf b)=Cl((\s)^{g/2})$. The case $g\equiv 0 \modu h_i$ has been treated in the first part of the proof, so we can suppose a fortiori that $\frac{g}{2}\not\equiv 0 \modu h_i$, thus the ideal $\mathbf b$ of $K$ is not principal. Therefore, the ideal $N_{K/L_2}(\s)$ of $L_2$ is not principal in the quadratic extension $L_2$ of $\Q$ between $\Q$ and $K$, and $h_i$ divides its class number $h(L_2)$, which achieves the proof.
\end{itemize}
\end{itemize}
\end{itemize}
\end{proof}
\end{thm}
%
{\bf Remark:} Observe that the class number $h(L_2)$  of the quadratic field $L_2/\Q, \quad \Q\subset L_2\subset K$ is small regards to the class number $h$ of $K/\Q$, which explains that in the numerical examples of the sequel, the $h_i$ candidate for the case $2.b)$ of the theorem are very small: in fact we have the upper bounds
\begin{itemize}
\item Quadratic real number fields:
  
$h(L_2)<\sqrt{D_2}$ where $D_2$ is the discriminant $D_2$ of the quadratic intermediate field $L_2$, see Narkiewicz, from proposition 8.7 p 451  \cite{nar}.
\item Quadratic imaginary fields:

We have the rough upper bound $h(L_2)<\sqrt{D_2}(Log(D_2))^2$ when $D_2>3000$ where $D_2$ is the discriminant of the quadratic intermediate field $L_2$, see Qu\^eme \cite{que}, proposition III.3 p 19-04.
\end{itemize}
%
\begin{thm}\label{t14051}
Let $K/\Q$ be a cyclic extension with $[K:\Q]=g$.
Let $g=2^\gamma\prod_{j=1}^m g_j^{\gamma_j}$ be the factorization of $g$ in primes.
Suppose that $d,\delta \in \N$, with $d\times\delta=g,\quad 2<d<g,\quad (d,\delta)=1$. 
Let $L_{d}$ be the intermediate field between $K$ and $\Q$ with $[L_{d}:\Q]=d$.
Let $D_{d}$ be the absolute discriminant of $L_{d}$.
Let $B_{d}$ be the Minkowski Bound of $L_d$.  
Let $G$ be the class group of $K/\Q$.
Let $h$ be the class number of $K/\Q$.
Let $h=\prod_{i=1}^r h_i^{\nu_i}$ be the factorization of $h$ in primes.
If for one prime factor $h_i$ of $h$, we have $h_i>(B_{d}\times(2Log(B_{d}))^d)$ and the $h_i$-component $G(h_i)$ of the class group $G$ of $K$  is cyclic, then  
$h_i\equiv 1 \modu g_j$ for a prime odd divisor $g_j$ of $\frac{g}{d}$.
\begin{proof}
The Galois group $Gal(K/L_d)$ is of order $\delta$. From $(d,\delta)=1$, we get $Gal(K/\Q)=G_1\oplus G_2$ where $G_1$ is a cyclic group of order $d$ and $G_2$ is a cyclic group of order $\delta$. We have $Gal(K/L_d)\cong G_2$ and therefore $Gal(K/L_d)$ is cyclic of order $\delta$.
Let $\tau:K \longrightarrow K$ be a $L_d$-isomorphism of $K$. Then the proof is similar to theorem (\ref{t22041}) proof, with the cyclic extension $K/L_d$ instead of the cyclic extension $K/\Q$. Here we have, for the ideal $\s$ of $K$ which generates a cyclic class of order $h_i$
\begin{displaymath}
\prod_{i=0}^{g/d-1} \tau^i (\s)=N_{K/L_d}(\s)B_d.
\end{displaymath}
The ideal $N_{K/L_d}(\s)B$ of the field $K$ generates a cyclic group of order $h_i$ or is principal.It cannot generate a cyclic group of order $h_i$ because $h_i$ cannot divide the order of the class of the ideal $N_{K/L_d}(\s)$ of $L_d$ : in fact  from hypothesis and lemma (\ref{l14051}) , $h_i$ is greater than the class number of $L_d$. The end of proof is then similar to theorem (\ref{t22041}) proof,  with $\frac{g}{d}$ in place of $g$. In this situation, the case $g\equiv 0 \modu h_i$ cannot occur, because roughly $(B_d(2Log(B_d))^d) >g$ and the case $h_i<\sqrt{D_2}$ cannot occur because roughly $(B_d\times(2Log(B_d))^d)>\sqrt{D_2}$.
\end{proof}
\end{thm}
%
\begin{thm}\label{t14052}
Let $K/\Q$ be a cyclic extension with $[K:\Q]=g$.
Let $g=2^\gamma\prod_{j=1}^m g_j^{\gamma_j}$ be the factorization of $g$ in primes.
Suppose that, for $j=1,\dots, m$,  we have $d_j=\frac{g}{g_j^{\gamma_j}}$. 
Let, for $j=1,\dots,m$,  $L_{d_j}$ be the intermediate field between $K$ and $\Q$ with $[L_{d_j}:\Q]=d_j$,
let $D_{d_j}$ be the absolute discriminant of $L_{d_j}$,
let $B_{d_j}$ be the Minkowski Bound of $L_{d_j}$.  
Let $G$ be the class group of $K/\Q$.
Let $h$ be the class number of $K/\Q$.
Let $h=\prod_{i=1}^r h_i^{\nu_i}$ be the factorization of $h$ in primes.
If for one prime factor $h_i$ of $h$, we have $h_i> Max(B_{d_j}\times(2Log(B_{d_j}))^{d_j}),\quad j=1,\dots m)$ and that the $h_i$-component $G(h_i)$ of the class group $G$ of $K$  is cyclic, then  
$h_i\equiv 1 \modu g$.
\begin{proof}
We apply theorem (\ref{t14051}) for $j=1,\dots,m$ observing that, for $j=1,\dots,m$,  we have $\frac{g}{d_j}=g_j^{\gamma_j}$, therefore $h_i\equiv 1 \modu g_j$ and then the result.
\end{proof}
\end{thm}
%
\begin{cor}\label{c22041}
Let $K/\Q$ be a cyclic extension with $[K:\Q]=g$ where $g\not\equiv 0 \modu 2$. 
Let $g=\prod_{j=1}^m g_j^{\gamma_j}$ be the factorization of $g$ in primes. Let $h$ be the class number of $K/\Q$.
Let $h=\prod_{i=1}^r h_i^{\nu_i}$ be the factorization of $h$ in primes.
If, for one prime factor $h_i$ of $h$, the $h_i$-component $G(h_i)$ of the class group $G$ of $K/\Q$ is cyclic, then we have  
\begin{itemize}
\item
else $g\equiv 0 \mod h_i$,
\item
else $g\not\equiv 0 \modu p$ and
$h_i\equiv 1 \modu g_j$ for a prime odd divisor $g_j$ of $g$, $1\leq j \leq m$.
\end{itemize}
\begin{proof} $ $
Immediate consequence of previous theorem (\ref{t22041}).
\end{proof}
\end{cor}
%
\begin{cor}\label{c22042}
Let $K/\Q$ be a cyclic extension with $[K:\Q]=g$ where $g$ is an odd prime. 
Let $h$ be the class number of $K/\Q$.
Let $h=\prod_{i=1}^r h_i^{\nu_i}$ be the factorization of $h$ in primes.
If for one prime factor $h_i$ of $h$, the $h_i$-component $G(h_i)$ of the class group $G$ of $K/\Q$ is cyclic, then we have 
\begin{itemize}
\item
else $h_i=g$,
\item
else $h_i\equiv 1 \modu g$.
\end{itemize}
\begin{proof} $ $
Immediate consequence of previous theorem (\ref{t22041}).
\end{proof}
\end{cor}
%
\begin{cor}\label{c22043}
Let $p\in \N$ be a prime such that $g=\frac{p-1}{2}$ is a prime.
Let $K=\Q(\xi_p+\xi_p^{-1})$ be the $p-$cyclotomic totally maximal real subfield of $\Q(\xi_p)$.
Let $h$ be the class number of $K/\Q$.
Let $h=\prod_{i=1}^r h_i^{\nu_i}$ be the factorization of $h$ in primes.
If for one prime factor $h_i$ of $h$, the $h_i$-component $G(h_i)$ of the class group $G$ of $K/\Q$ is cyclic, then we have
\begin{itemize}
\item
else $h_i=g$,
\item
else $h_i\equiv 1 \modu g$.
\end{itemize}
\begin{proof} $ $
Immediate consequence of previous theorem (\ref{t22041}).
\end{proof}
\end{cor}

{\bf Remark: } Note that when $h_i\| h$, the condition $G(h_i)$ is a cyclic subgroup of $G$ is always verified.
%
%
\subsection{Numerical examples}\label{ss22043}
The examples found to  check theses results are taken from: 
\begin{itemize} 
\item
the table of relative class numbers of cyclotomic number fields in Whashington, \cite{was} p 412,  with some elementary MAPLE computations,
\item 
the table of relative class number of cyclotomic number fields in School, \cite{sch}
\item
the table of maximal real subfields $\Q(\xi_p+\xi_p^{-1})$ of $\Q(\xi_p)$ for $p$ prime in Washington, \cite{was} p 420.
\item 
the tables of relative class number of imaginary cyclic fields of Girstmair of degree 4,6,8,10. \cite{gir}.
\item
the tables of quintic number fields computed by Jeannin, \cite{jea}
\item 
the tables of cubic totally real cyclic number fields of the Bordeaux University in the Server 

{\it megrez.math.u-bordeaux.fr}.
\end{itemize}

All the results examined in these tables are in accordance with our theorems.
\subsubsection{Cyclotomic number fields $\Q(\xi_n)$}

The cyclotomic number fields $\Q(\xi_n)$ of the example are  with $2\| \phi(n)$.
Here we have $g=\phi(n)$. 
\begin{itemize}
\item 

$\Q(\xi_n), \quad n=59 ,\quad g=\phi(n)=58=2.29 :$

$h^- = 3 . (2.29+1). (2^3 . 29 +1)$

$3 =h(\Q(\sqrt{-59})$
\item

$\Q(\xi_n), \quad n=71, \quad g=\phi(n)=70=2.5.7$ :

$h^- = 7^2(2^3.5.7.283+1)$.

\item
$\Q(\xi_n), \quad n=79,\quad  g=\phi(n)=78=2.3.13$ :

$ h^-= 5.(2^2.13+1)(2.3^2.5.13.17.19+1)$.

$5=h(\Q(\sqrt{-79})$.

\item
$\Q(\xi_n), \quad n=83,\quad g=\phi(n)=82=2.41$ :

$ h^-=3.(2^2.41.1703693+1)$.

$3=h(\Q(\sqrt{-83}).$
\item
$\Q(\xi_n),\quad n=103,\quad  g=\phi(n)=102=2.3.17$:

$ h^-= 5.(2.3.17+1)(2^2.3.5.17+1)(2.3^2.5.17.11273+1)$.

$5=h(\Q(\sqrt{-103})$
\item
$\Q(\xi_n), \quad n=107,\quad g=\phi(n)=106=2.53$ :

$ h^-=(2.7.53+1)(2.3.31.53+1)(2^6.23.37.53+1)$.
\item
$\Q(\xi_n), \quad n=121,\quad g=\phi(n)=110=2.5.11$ :

$ h^-=(2.3.11+1)(2^5.11+1)(2^2.5.11.13+1)(2^2.3^2.5.11.13)$.

\item
$\Q(\xi_n), \quad n=127,\quad g=\phi(n)=126=2.3^2.7$ :

$ h^-=5.(2^2.3+1)(2.3.7+1)(2.3.7.13+1)(2.3^2.7^2+1)\times$

$5=h(\Q(\sqrt{127})$.

$(2.3^4.19+1)(2.3^2.7.4973+1)$.
\item
$\Q(\xi_n), \quad n=131,\quad g=\phi(n)=130=2.5.13$ :

$ h^-=3^3.5^2.(2^2.13+1)(2.5.13+1)(2^2.5^2.13+1)\times$

$(2^2.5^2.13.29.151.821+1)$.

For $h_i=3$, we have  $h_i^3\|h$ and the group $G(h_i)$ is not cyclic as it is seen in School \cite {sch} table 4.2 p 1239 : we cannot apply theorem (\ref{t22041}) to this number $h_i=3$.

\item
$\Q(\xi_n), \quad n=139,\quad g=\phi(n)=138=2.3.23$ :

$ h^-=3^2(2.23+1)(2^2.3.23+1)(2.3.7.23+1)(2^2.3.23.4307833+1)$.
\item
$\Q(\xi_n), \quad n=151,\quad g=\phi(n)=150=2.3.5^2$ :

$ h^-=(2.3+1)(2.5+1)^2.(2^2.5.7+1)(2.3.5^2.173+1)\times$

Observe that for $h_i=2.5+1$, we cannot conclude because we have not proved that $G(h_i)$ is cyclic. Moreover, from Schoof, see table 4.2 p 1239, that the group is effectively not cyclic.

$(2^2.3.5^4.7.23+1)(2^2.3.5^2.7.13.73.1571+1)$.
\item
$\Q(\xi_n), \quad n=163,\quad g=\phi(n)=150=2.3^4$ :

$ h^-=(2^2.3^2.5+1)(2.3^4.11.13+1)(2^5.3^5.47+1)\times$

$(2.3^4.17.19.29.71.73.56179+1)$.
\item
$\Q(\xi_n), \quad n=167,\quad g=\phi(n)=166=2.83$ :

$ h^-=11.(2.3.83+1)(2.83.22107011.1396054413416693+1)$.

$11=h(\Q(\sqrt{167})$.
\item
$\Q(\xi_n), \quad n=179,\quad g=\phi(n)=178=2.89$ :

$ h^-=(2^2.3.89+1)\times$

$(2^4.5.89.173.19207.155731.3924348446411+1)$.

\item $\Q(\xi_n), \quad n=191,\quad g=\phi(n)=190=2.5.19$:

$ h^-= (2.5+1).13.(2.19^2.71+1)(2.3.5.19.277.3881+1)$.

$13=h(\Q(\sqrt{-191})$.
\item

$\Q(\xi_n), \quad n=199,\quad g=\phi(n)=198=2.3^2.11$ :

$ h^-=3^4(2.3^2+1)(2.3.11^2+1)(2^2.3.11.23.8447)(2^4.3^2.11.13.17.331.1789)$.
\item $\Q(\xi_n), \quad n=211, \quad g=\phi(n)=2.3.5.7$,
  
$h^- = 3 . (2.3+1)(2^2.5+1)(2^2. 3^2 5+1)(2^3. 5.7+1)
(2 .3.5^2.7+1)\times$

$(2.5^3.7+1)(2^2.3.5^2.7.54277+1)
(2^2.3.5.59.122743+1)$
\item $Q(\xi_n),\quad n=223,\quad g=\phi(n)=222=2.3.37$:

$h^-=(2.3+1)(2.3.7+1)(2.13.17.37.34141.32077+1)\times$

$(2.3.5^2.37.419417.5051)$.
\item $Q(\xi_n),\quad n=227,\quad g=\phi(n)=226=2.113$:

$h^-=5.(2.13.113+1)^3.(2^2.5^2.113.149807435573+1)\times$

$(2^2.3.113.2207.4973.903334373+1)$.

$5=h(\sqrt{\Q(-227})$.
\item
$\Q(\xi_n),\quad n=239,\quad g=\phi(n)=238=2.7.17$:

$ h^- = $

$3.5.(2.3.17.5011+1)(2.3.13.17.1523+1)\times$

$(2^2.7.17.41.6321643+1)\times$

$(2.3.7.17.31.182453.315075469.17681+1)$.

$3.5=h(\Q(\sqrt{239})$.

\item $\Q(\xi_n), \quad n=251,\quad  g=\phi(n)=250=2.5^3$:

$h^-= 7.(2.5+1)(2^4.3.5^4.13.19.47+1)\times$

$(2.3^2.5^3.89.1559.30851+1)\times$

$(2^2.3^2.5^3.463.827.214520849205730542617+1).$

$7=h(\Q(\sqrt{251})$.
\item
$\Q(\xi_n), \quad n=307,\quad  g=\phi(n)=306=2.3^2.17$:

$h^-= 3^3(2^2.3^2+1)(2^3.17+1)(2.13.17+1)\times$

$(2^2.3.7.17+1)(2.3^2.17+1)(2^5.3^2.13.17.160119221+1)\times$

$(2^2.3^2.17+1)(2.3^3.17+1)\times$

$(2.3^3.17.37459802253317.11580768409391.1286690597+1)$.

\end{itemize}
{\bf Remarks:} 
\begin{itemize}
\item
We see that almost all these results correspond to the first and second case of theorem (\ref{t22041}) and few to the third case.  
\item
We see also that the hypothesis $G(h_i)$ cyclic is verified for a great part of the primes $h_i$ and all $h_i$ sufficiently large.
\item 
We see that the results of theorems (\ref{t14051}) and (\ref{t14052}) are well verified.
\end{itemize}
%
\subsection{Real class number $\Q(\xi_p+\xi^{-1}_p)$}
The examples are obtained from the table of real class number in Washington, \cite{was} p 420.
Here, $p$ is a prime, the class number  $h_\delta$ is the conjectured value of the class number $h^+$ of 
$\Q(\xi_p+\xi_p^{-1})/\Q$ with a minor incertitude on an extra factor. But the factor $h_\delta$ must verify our theorems. 
We extract some examples of the table with $2\| p-1$. Here we have $g=\frac{p-1}{2}$ and thus $g\not \equiv 0 \modu 2$. Therefore the factor $h_i=2$ cannot be taken in account  in the corollary (\ref{c22041}) which implies that $G(2)$ is not cyclic and that $Card(G(2))=2^\alpha, \alpha\geq 2$, which is verified in the Table.
\begin{itemize}
\item
$p=163,\quad p-1=2.3^4,\quad h_\delta=2^2, \quad G(2)$ not cyclic.

\item
$p=191,\quad p-1=2.5.19,\quad h_\delta=(2.5+1)$

\item
$p=491,\quad p-1=2.5.7^2,\quad h_\delta=2^3, \quad G(2)$ not cyclic.

\item
$p=547,\quad p-1=2.3.5.13,\quad h_\delta=2^2, \quad G(2)$ not cyclic.

\item
$p=827,\quad p-1=2.7.59,\quad h_\delta=2^3, \quad G(2)$ not cyclic.

\item
$p=1063,\quad p-1=2.3^2.59,\quad h_\delta=(2.3+1).$

\item
$p=1231,\quad p-1=2.3.5.41,\quad h_\delta=(2.3.5.7+1).$

\item
$p=1399,\quad p-1=2.3.233,\quad h_\delta=2^2, \quad G(2)$ not cyclic.

\item
$p=1459,\quad p-1=2.3^6,\quad h_\delta=(2.3.41+1).$

\item
$p=1567,\quad p-1=2.3^3.29,\quad h_\delta=(2.3+1).$

\item
$p=2659,\quad p-1=2.3.443,\quad h_\delta=(2.3^2+1).$

\item
$p=3547,\quad p-1=2.3^2.197,\quad h_\delta=((2.3^2+1)(2.3^2.7^2+1).$

\item
$p=8017,\quad p-1=2^4.3.167$,

$h_\delta=(2.3^2+1)(2.3^2.7^2+1)(2^2.3^3+1).$

\item
$p=8563,\quad p-1=2.3.1427,\quad h_\delta=(2.3+1)^2$, we cannot conclude because we cannot prove that $G(7)$ is cyclic.

\item
$p=9907,\quad p-1=2.3.13.127,\quad h_\delta=(2.3.5+1).$
\end{itemize}

%

\subsubsection {Cubic fields $K/\Q$ cyclic and totally real.}
Here, we have $g=3$. Note that in that case discriminants are square in $\N$. 
\begin{itemize}
\item Discriminant $D_1=3969=(3^2.7)^2$, \quad$h=3$
\item Discriminant $D_2=3969=(3^2.7)^2$, \quad $h=7=2.3+1$
\item Discriminant $D_1=8281=(7.13)^2$, \quad $h=3$
\item Discriminant $D_2=8281=(7.13)^2$, \quad $h=3$
\item Discriminant $D_1=13689=(3^2.13)^2$,\quad $h=3$
\item Discriminant $D_2=13689=(3^2.13)^2$, \quad $h=13=2^2.3+1$
\item Discriminant $D_1=17689=(7.19)^2$, \quad $h=3$
\item Discriminant $D_2=17689=(7.19)^2$, \quad $h=3$
\item Discriminant $D=26569=163^2$, \quad $h=4=2^2$, we cannot apply theorem because the class group
$C_2\times C_2$ is not cyclic.
\item Discriminant $D_1=47089=(7.31)^2$,\quad  $h=3$
\item Discriminant $D_2=47089=(7.31)^2$,\quad  $h=3$
\item Discriminant $D_1=61009=(13.19)^2$,\quad  $h=3$
\item Discriminant $D_2=61009=(13.19)^2$, \quad $h=3$
\item Discriminant $D_2=67081=(7.37)^2$\quad  $h=3$
\item Discriminant $D_1=67081=(7.37)^2$\quad  $h=3$
\item Discriminant $D_1=76729=277^2$\quad  $h=2\times 2$

$G(h_1)$ not cyclic
\item Discriminant $D_1=77841=(3^2.31)^2$\quad  $h=3$
\item Discriminant $D_2=77841=(3^2.31)^2$\quad  $h=3$
\item Discriminant $D_1=90601=(7.43)^2$, \quad $h=3$
\item Discriminant $D_2=90601=(7.43)^2$, \quad $h=3$
\item Discriminant $D=97969=313^2$, \quad $h=7=2.3+1$
\item Discriminant $D_1=110889=(3^2.37)^2$,\quad  $h=3$
\item Discriminant $D_2=110889=(3^2.37)^2$,\quad $h=3$
\item Discriminant $D=121801=349^2$, \quad $h=4=2^2$

$G(h_1)$ not cyclic

\item Discriminant $D_1=149769=(3^2.43)^2$,\quad  $h=3$
\item Discriminant $D_2=149769=(3^2.43)^2$, $h=3$
\item Discriminant $D=157609=397^2$, \quad $h=4=2^2$

$G(h_1)$ not cyclic
\item Discriminant $D_1=162409=(13.31)^2$,\quad  $h=3$
\item Discriminant $D_2=162409=(13.31)^2$, \quad $h=3$
\item Discriminant $D_1=167081=(7.37)^2$\quad  $h=3$
\item Discriminant $D_1=182329=(7.61)^2$\quad  $h=3$
\item Discriminant $D_2=182329=(7.61)^2$\quad  $h=3$
\item Discriminant $D_1=368449=607^2$\quad  $h=2\times 2$

$G(h_1)$ not cyclic
\item Discriminant $D_1=727609=853^2$\quad  $h=2\times 2$

$G(h_1)$ not cyclic
\item Discriminant $D_1=1482809=(3^2.7.19)^2$\quad  $h=3\times 3$

$G(h_1)$ not cyclic

\item Discriminant $D_2=1482809=(3^2.7.19)^2$\quad  $h=6\times 6$

$G(h_1), G(h_2)$ not cyclic
\item Discriminant $D_3=1482809=(3^2.7.19)^2$\quad  $h=3\times 3$

$G(h_1)$ not cyclic

\item Discriminant $D_4=1482809=(3^2.7.19)^2$\quad  $h=3\times 3$

$G(h_1)$ not cyclic
\end{itemize}
%
\subsection{Totally real cyclic fields of prime conductor $<100$}
We have found few numeric results in the literature. We refer to Masley, \cite{mas}, Table 3 p 316.
In these results, $K/\Q$ is a real cyclic field with $[K:\Q]=g$, with conductor $f$, with root of discrimant $Rd$ and class number $h$.
\begin{itemize}
\item 
$f=40,\quad g=2,\quad \quad h=2$
\item 
$f=60,\quad g=2,\quad \quad h=2$
\item 
$f=63,\quad g=3,\quad Rd=15.84,\quad \quad h=3$
\item 
$f=63,\quad g=3,\quad Rd=15.84,\quad \quad h=3$
\item 
$f=63,\quad g=6,\quad Rd=26.30,\quad \quad h=3$
\item 
$f=63,\quad g=6,\quad Rd=26.30,\quad \quad h=3$
\item 
$f=65,\quad g=2,\quad \quad h=2$
\item 
$f=91,\quad g=3,\quad Rd=20.24,\quad \quad h=3$
\item 
$f=91,\quad g=3,\quad Rd=20.24,\quad \quad h=3$
\item 
$f=91,\quad g=6,\quad Rd=31.03,\quad \quad h=3$
\item 
$f=,\quad g=6,\quad Rd=31.03,\quad \quad h=3$
\end{itemize}
%
\subsection{Lehmer quintic cyclic field}
The prime divisors of the $82$ cyclic number fields of the table in Jeannnin \cite{jea}, with conductor $f<3000000$,  are, at a glance, of the form $h_i=2$ or $h_i=5$ or $h_i=2.5.k+1$, which clearly verifies theorem (\ref{t22041}). The corollary (\ref{c22041}) cannot be applied to the prime $h_i=2$ because, here $g=5\not \equiv 0 \modu 2$ : therefore, we deduce that the $2-$component of the class-group $G(2)$ cannot be cyclic : we observe that the power of $2$ in the class number is $2^4$ or $2^8$ which is compatible with $G(2)$ not cyclic. 
%
\subsection{Decimic imaginary cyclic number fields with conductor between  $9000$ and $9500$}
This example is obtained from the tables of Girstmair, \cite{gir}. $p$ is a prime conductor, $h$ is the factorization of the class number $K/\Q$.
\begin{itemize}
\item 
$f=9011$, $h=3.11.1566031$.

$3$ divides the class number of $Q(\sqrt{-9011})$.
\item 
$f=9081$, $h=3.7.41.491$.

$3,7$ divide the class number of $Q(\sqrt{-9081})$.
\item 
$f=9151$, $h=67.23741$.

$67$ divides the class number of $Q(\sqrt{-9151})$.
\item 
$f=9311$, $h=97.56891$.

$97$ divides the class number of $Q(\sqrt{-9311})$.
\item 
$f=9371$, $h=7^2.151.271$.

$7$ divides the class number of $Q(\sqrt{-9371})$.
\item 
$f=9391$, $h=5^2.119281$.

\item 
$f=9431$, $h=7.13.31.17041$.

$7,13$ divide the class number of $Q(\sqrt{-9431})$.
\item 
$f=9491$, $h=3^2.5^2.31.2341$.

$3$ divides the class number of $Q(\sqrt{-9491})$.
\end{itemize} 
%
\subsection{A commentary on results}
We use the hypothesis $G(h_i)$ is a cyclic group. 
\begin{itemize}
\item
We observe that for several $G(h_i)$ not cyclic we have nevertheless $h_i\equiv 1 \modu g_j$.
An example very significative is 
$\Q(\xi_p)$ with $p=227$. We have $p-1=2.113$. Then $g_1=113$. The group 
$G(2939)=C_{2939}\times C_{2939}\times C_{2939}$ is not cyclic, see Schoof \cite{sch} table 4.2 p 1239. But we have $h_i=2939=(2.13.113+1)$.
\item An example a contrario $\Q(\xi_p)$ with $p=131$. here $p-1=2.5.13$.
$G(3)=C_3\times C_3 \times C_3$ is not cyclic (see Schoof), moreover $3$ does not divides $h(\sqrt{-263})=5$. The theorem (\ref{t22041}) is not applicable. We have $h_i^{\nu_i}=27=(2.13+1)$.
here $\nu_i=3=$ rank of $G(3)$.
\item Another example a contrario. Let $p=263, \quad p-1=2.131$. From Schoof \cite{sch} p1243, $3^4\| h^-$. We can assert that the group $G(3^4)$ is not cyclic, because from theorem(\ref{t22041}), we would get $3 | h(\Q(\sqrt{-263}=13$. In an other part we see that the smaller $\nu$ such that $3^\nu\equiv 1 \modu 131$ is $\nu=65$. This $\nu$ is greater than the rank of $G(3)$ which is $2$ or $4$.
\end{itemize} 
%
\section{On class number of cyclotomic extensions}\label{s22042}
This section gives another elementary method of study of factorization of class number of cyclotomics number fields. The approach is strictly independant of the previous part of this article.
\subsection{Some definitions and notations}\label{s1}
\begin{itemize}
\item 
Let $p\in \N,\quad p>3$ be a  prime.
\item
Let $\xi_p$ with $\xi_p^{p-1}+\xi_p^{p-2}+\dots+\xi_p+1=0$, where $\xi_p\in\C$, is a primitive root.
\item
Let $\Q(\xi_p)$ be the $p-$cyclotomic number field.
\item 
Let $\Z[\xi_p]$ be the ring of integers of $\Q(\xi_p)$.
\item
Let $\Z[\xi_p]^*$ be the group of units of $\Z[\xi_p]$.
\item 
Let $a,b\in\Z-\{0\},\quad (a+b)\not=0,\quad \frac{(a^p+b^p)}{(a+b)}>1$.
Let us suppose that $a,b$ are  coprime with $ab(a+b)\not\equiv 0 \modu p$.
\item
Let $n=p\frac{(a^p+b^p)}{(a+b)}$.
\item
Let $\Q(\xi_n)$ be the $n-$cyclotomic extension of $\Q$.
\item 
Let $n=p\prod_{i=1}^f q_i^{\alpha_i}$ be the factorization in primes in $\N$. Classically, we have $q_i\equiv 1 \modu p,\quad i=1,\dots,f$.
\item
We have $\phi(n)=(p-1)\prod_{i=1}^f \phi(q_i^{\alpha_i})$.
\item
Note that $2$ is prime with $n$, so $(n,2)=1$.
\item 
Let $\Z[\xi_n]$ be the ring of integers of $\Q(\xi_n)$.
\item
Let $\Z[\xi_n]^*$ be the group of units of $\Q(\xi_n)$.
\item
Let $\Q(\xi_n+\xi_n^{-1})$ be the totally real  maximal subfied of $\Q(\xi_n)$.
\item
Let $\Z[\xi_n+\xi_n]$ be the ring of integers of $\Q(\xi_n)$.
\item
Let $\Z[\xi_n+\xi_n]^*$ be the group of units of $\Z[\xi_n+\xi_n^{-1}]$.
\end{itemize}

\subsection{On class number of the cyclotomic $n-$field $\Q(\xi_n)$.}
\begin{lem}
Let $a,b\in\Z-\{0\},\quad (a+b)\not=0,\quad\frac{(a^p+b^p)}{(a+b)}>0$. 
Let us suppose that $a,b$ are  coprime with $ab(a+b)\not\equiv0 \modu p$.
Let $q,\quad q\in\N$ be a prime dividing $t_1=\frac{(a^p+b^p)}{(a+b)}$. 
Let $\mathbf q_p$ be a prime ideal of $\Q(\xi_p)$ above $q$. Then $q\equiv 1 \modu p$
and $\mathbf q_p\Z[\xi_{pq}]=\mathbf q_{p,q}^{q-1}$ where $\mathbf q_{pq}$ is a prime ideal of $\Z[\xi_{pq}]$. 
\begin{proof}
$ $
\begin{itemize}
\item 
From hypothesis on $a,b$, we get $t_1\not\in\{-1,1\}$. Therefore, there exists at least one prime $q\in\N$ dividing $t_1$.
\item
From $t_1\equiv 0 \modu q$, we get classically $q\equiv 1 \modu p$.
\item
Decomposition of $q$ in $\Q(\xi_q)$  : There is one  prime $\mathbf q_{q}=(1-\xi_q)\Z[\xi_q]$ above $q$ in $\Q(\xi_q)$ with inertial degree 1, ramification indice $q-1$.
\item 
Decomposition of $q$ in $\Q(\xi_p)$ : there are $p-1$ ideals $q_{p,i},\quad
i=1,\dots,p-1$ with inertial degree $1$ and ramification indice $1$.
\item 
Decomposition of $q$ in $\Q(\xi_{pq})$ : 
\begin{itemize}
\item
The ramification indices and inertial degrees are multiplicative in tower and $[\Q(\xi_{pq}):\Q]=(p-1)(q-1)$;
\item
therefore there is only on possible configuration : there are $p-1$ ideals $\mathbf q_{pq,i},\quad i=1,\dots,p-1$ of inertial degree $1$ and ramification indice $q-1$ 
above $q_{p,i}$.
\end{itemize}
\item
Then $\mathbf q_{p,i}\Z[\xi_{pq}]=\mathbf q_{pq,i}^{q-1}$ where $q-1\equiv 0\modu p$.
\end{itemize}
\end{proof}
\end{lem}  
\begin{thm}
Let $a,b\in\Z-\{0\},\quad (a+b)\not=0,\quad \frac{(a^p+b^p)}{(a+b)}>1$. 
Let us suppose that $a,b$ are co-prime with $ab(a+b)\not\equiv 0 \modu p$. 
Let $n=p\frac{(a^p+b^p)}{(a+b)}$.
 Then
$(a+\xi_p b)\Z[\xi_n]=\mathbf U^p$, where $\mathbf U$ is an integral ideal of $\Z[\xi_n]$. 
\begin{proof}
$ $
\begin{itemize}
\item
Let $\mathbf q_p$ be a prime ideal of $\Z[\xi_p]$ above $q$. From previous lemma,
$\mathbf q_p\Z[\xi_{pq}]=\mathbf q_{pq}^{q-1}$, where $\mathbf q_{pq}$ is a prime ideal of $\Z[\xi_{pq}]$. Then, from $q\equiv 1 \modu p$ we get 
$\mathbf q_p\Z[\xi_{pq}]=\mathbf v_q^p$ where $\mathbf v_q$ is an integral ideal of $\Z[\xi_{pq}]$.
\item
We can extend this result to all prime $q\not=p$ which divides $n$. We note that $\Z[\xi_{pq}]\subset \Z[\xi_n]$ for all prime $q\not=p$ dividing $n$, which then leads to the result.
\end{itemize}
\end{proof}
\end{thm}
\begin{thm}\label{t28115}
Let $a,b\in\Z-\{0\},\quad (a+b)\not=0, \quad\frac{(a^p+b^p)}{(a+b)}>1$.
Let us suppose that $a,b$ are  coprime with $ab(a+b)\not\equiv 0\modu p$. 
Let $n=p\frac{(a^p+b^p)}{(a+b)}$. Let $\Q(\xi_n)$ be the $n-$cyclotomic field.
Let $h_n$ be the class number of $\Q(\xi_n)$.
Then $h_n\equiv 0 \modu p$.
\begin{proof}
We shall suppose that $h_n\not\equiv 0\modu p$ to obtain a contradiction.
\begin{itemize}
\item
From previous theorem, we have $(a+\xi_p b)\Z[\xi_n]=\mathbf U^p$, where $\mathbf U$ is an integral ideal of $\Z[\xi_n]$. From $h_n\not\equiv 0 \modu p$, then  $\mathbf U$ is a principal integral ideal of $\Z[\xi_n]$.
Then we can write $(a+\xi_p b)=\varepsilon \gamma^p$ where $\varepsilon\in\Z[\xi_n]^*$ and where $\gamma\in\Z[\xi_n]$.
\item
The cyclotomic $n-$field $\Q(\xi_n)/\Q$ is a complex normal extension. From a generalization of Kummer Lemma, see \cite{nar} proposition 3.6 p 109, we have  $\varepsilon=\theta\eta$ where $\theta^2\in \Q(\xi_n),\quad \theta^2=\xi_n^{l_0}$ and where $\eta^2\in \Z(\xi_n+\xi_n^{-1})^*$.
Therefore,
\begin{equation}
(a+\xi_p b)^2= \xi_n^{l_0}\eta_1\gamma_1^{2p},\quad \eta_1 \in\Z[\xi_n+\xi_n^{-1}]^*, \quad \gamma_1\in\Z[\xi_n]
\end{equation}
\item
By conjugation, we obtain
\begin{equation}
(a+\xi_p^{-1}b)^2=\xi_n^{-l_0}\eta_1\overline{\gamma_1}^{2p}
\end{equation}
Then 
\begin{equation}
\frac{(a+\xi_p b)^2}{(a+\xi_p^{-1}b)^2} =\xi_n^{2l_0}(\frac{\gamma_1}{\overline{\gamma_1}})^{2p}.
\end{equation}
which can be written
\begin{equation}
\frac{(a+\xi_p b)}{(a+\xi_p^{-1} b)}=\xi_n^{l_1} C_1^p, \quad C_1\in\Q(\xi_n),
\quad 0 \leq l_1 \leq n.
\end{equation}
\item
Let us define the map $\xi_n\Longrightarrow\xi_n^2\in Gal(\Q(\xi_n)/\Q).$
\begin{equation}
\frac{(a+\xi_p^{2} b)}{(a+\xi_p^{-2} b)}=\xi_n^{2 l_1} C_2^p, \quad C_2\in\Q(\xi_n),
\quad 0 \leq l_2 \leq n.
\end{equation}
\item 
From this relation, by division, we get
\begin{displaymath}
\frac{(a+\xi_p b)^{2} (a+\xi_p^{-2} b)}
     {(a+\xi_p^{-1}b)^{2} (a+\xi_p^{2} b)}= (\frac{C_1^{2}}{C_2})^p
\end{displaymath}
Therefore, there exists $C\in\Q(\xi_n)$ such that
\begin{displaymath}
\frac{(a+\xi_p b)^{2} (a+\xi_p^{-2} b)}
     {(a+\xi_p^{-1}b)^{2} (a+\xi_p^{2} b)}= C^p.
\end{displaymath}
Therefore, at this stage of the proof, we have simultaneously $C\in\Q(\xi_n)$ and $C^p\in\Q(\xi_p)$.
\item 
We have $n=p\prod_{i=1}^f q_i^{\alpha_i}$ where the prime $q_i,\quad i=1,\dots,f$ verify 
$q_i\equiv 1 \modu p$. Then $\phi(n)=\phi(p)\times \prod_{i=1}^f \phi(q_i)^{\alpha_i}$ leads to
$\phi(n)=p^\alpha\times m,\quad m\not\equiv 0 \modu p, \quad \alpha\geq 1$.
\item There is one and only one Kummer extension $K/\Q(\xi_p)$ with $[K:\Q(\xi_p)]=p$ and with
$K\subset \Q(\xi_n)$, given by $K=\Q(\xi_p,\xi_{p^2})=\Q(\xi_{p^2})$. 
\item There are two possibilities $C\in\Q(\xi_p)$ or $C\in\Q(\xi_n)-\Q(\xi_p)$. If $C\in\Q(\xi_n)-\Q(\xi_p)$ then, from $C^p\in\Q(\xi_p)$, we deduce that
$\Q(\xi_p,C)/\Q(\xi_p)$ is a Kummer extension with $[\Q(\xi_p,C):\Q(\xi_p)]=p$. Therefore, from $C\in\Q(\xi_n)$, we get
$Q(\xi_p,C)=\Q(\xi_{p^2})$ and hence $C=\xi_{p^2}^u C_u$ where $u\in\Z,\quad C_u\in\Q(\xi_p)$.
\item 
 Therefore
we have $u\in\Z,\quad C_u\in\Q(\xi_p)$ with
\begin{displaymath}
\frac{(a+\xi_p b)^{2} (a+\xi_p^{-2} b)}
     {(a+\xi_p^{-1}b)^{2} (a+\xi_p^{2} b)}= (\xi_{p^2}^u C_u)^p=\xi_p^u C_u^p.
\end{displaymath}
\item
Then, we have 
\begin{displaymath}
\frac{(a+\xi_p^{-1} b)^{2} (a+\xi_p^{2} b)}
     {(a+\xi_p^b)^{2} (a+\xi_p^{-2} b)}= \xi_p^{-u}\overline{C_u}^p.
\end{displaymath}
Classically $C_u^p-\overline{C_u}^p\equiv 0 \modu p$. 
Then
\begin{displaymath}
\xi_p^{-u}\frac{(a+\xi_p b)^{2} (a+\xi_p^{-2} b)}
     {(a+\xi_p^{-1}b)^{2} (a+\xi_p^{2} b)}\equiv
\xi_p^u\frac{(a+\xi_p^{-1} b)^{2} (a+\xi_p^{2} b)}
     {(a+\xi_p b)^{2} (a+\xi_p^{-2} b)} \modu p,
\end{displaymath}
hence
\begin{displaymath}
\xi_p^{-u}(a+\xi_p b)^{2} (a+\xi_p^{-2} b)
     \equiv
\xi_p^u(a+\xi_p^{-1} b)^{2} (a+\xi_p^{2} b)
     \modu p.
\end{displaymath}
\item
Then, with hypothesis $p>3$,  we conclude, clearly, that this congruence is not possible.
\item
Then $h_n\not\equiv 0$ would lead to a contradiction, which achieves the proof.
\end{itemize}
\end{proof}
\end{thm}
\subsection{Numerical examples}\label{ss22043}
In this section, we give examples in accordance with Washington \cite{was} p 412 table of relative class number of the cyclotomic fields $\Q(\xi_n)$, for all $n,\quad \phi(n)\leq 256$.
\begin{itemize}
\item 
$p=5,\quad a=2,\quad b=1,\quad n=55,\quad \phi(n)= 40,\quad h_n^-=2\times 5 \equiv 0 \modu 5$.
\item 
$p=5,\quad a=3,\quad b=1,\quad n=305,\quad \phi(n)= 240,\quad h_n^-=2^3\times 3^2\times 
5^4\times 13^2\times \dots \equiv 0 \modu 5$.
\item 
$p=7,\quad a=2,\quad b=1,\quad n=301,\quad \phi(n)= 252,\quad h_n^-=2^{10}\times
3^3\times 7^7\times 19\times \dots  \equiv 0 \modu 7$.
\end{itemize}
\section{Comparison with bibliography}
\begin{itemize}
\item
We have not identified the theorem (\ref{t22041}) and corollaries (\ref{c22041}), (\ref{c22042}) and (\ref{c22043}) explicitly in this direct form in the bibliography. In relation with  these results, the next proposition in Gras, \cite {gra} :
Let $L/\Q$ be a cyclic extension.
Let $\chi$ be a rational character, $\chi\in \Xi_L$, where the set $\Xi_L$ is defined in Gras \cite{gra} p 36. Let, as usual, the symbol $+$ and $-$ for the real part, and relative part of $L/\Q$. Let $H(L),H(L)^+, H(L)^-$ be the set of class of $L/\Q$, respectively $L^+/\Q, L^-/\Q$. Let $H_\chi$ indexed by $\chi\in \Xi_L$, $H^\prime_\chi$ with the definition of the correspondance $O \Longrightarrow O^\prime$ between objects $O$ and $O^\prime$ defined in  \cite{gra}, p 36. Then we have, in \cite{gra}, proposition II.2 p 39 : 
\begin{displaymath}
\begin{split}
&|H(L)|=\prod_{\chi\in \Xi_L} |H_\chi^\prime|,\\
&|H(L)^+|=\prod_{\chi\in \Xi_L^+} |H_\chi^\prime|,\\
&|H(L)^-|=\prod_{\chi\in \Xi_L^-} |H_\chi|,\\
\end{split}
\end{displaymath}
It is noteworthy that our approach is strictly elementary.
\item
The theorem (\ref{t28115}) is a weaker form  of results on genus theory, for instance, it is possible to prove the stronger result that if $p$ and $q$ are two prime numbers and $h_{pq}$ is the class number of  $\Q(\xi_{pq})$ and  if $q\equiv 1 \modu p$ then $h_{pq}\equiv 0 \modu p$ in the two cases: $p=3$ and $q\equiv 1 \modu 9$ and $p\geq 5$. In  this direction see Ishida \cite{ish}, see G. Gras \cite{gr2}. By opposite,the  proof of theorem (\ref{t28115}) is strictly elementary.
\end{itemize}

%

Roland Qu\^eme

13 avenue du ch\^ateau d'eau

31490 BRAX

FRANCE

e-mail : 106104.1447@compuserve.com

\end{document}